\documentclass[12pt,a4paper]{article}
\usepackage[english]{babel}
\usepackage{amssymb}
\usepackage[T1]{fontenc}
\usepackage{amsmath}
\newtheorem{Th}{Theorem}
\newtheorem{lemma}{Lemma}
\newtheorem{Cor}{Corollary}
\newtheorem{Prop}{Proposition}
\newtheorem{RMQ}{Remark}

\newcommand{\F}{\mathcal{F}}
\newcommand{\R}{\mathbb{R}}
\newcommand{\Z}{\mathbb{Z}}
\newcommand{\N}{\mathbb{N}}

\newcommand{\U}{\mathcal{U}}
\newcommand{\V}{\mathcal{V}}

\renewcommand{\P}{\mathbb{P}}

\newcounter{tictac}
\newenvironment{fleuveA}{
   \begin{list}{$\textbf{A\arabic{tictac}}$) }{\usecounter{tictac}
\leftmargin 1cm\labelwidth 2em}}{\end{list}}
\newenvironment{fleuve}{
   \begin{list}{$\textbf{\emph{\arabic{tictac})}}$ }{\usecounter{tictac} \leftmargin 1cm\labelwidth
2em}}{\end{list}}
\newenvironment{fleuveC}{
   \begin{list}{$\textbf{C\arabic{tictac}}$) }{\usecounter{tictac}
\leftmargin 1cm\labelwidth 2em}}{\end{list}}
\def\1{\,\rlap{\mbox{\small\rm 1}}\kern.15em 1}
\def\ind#1{\1_{#1}}
\def\build#1_#2^#3{\mathrel{\mathop{\kern 0pt#1}\limits_{#2}^{#3}}}
\def\tend#1#2{\build\hbox to 12mm{\rightarrowfill}_{#1\rightarrow #2}^{a.s.}}

\def\converge#1#2#3{\build\hbox to
15mm{\rightarrowfill}_{#1\rightarrow #2}^{\hbox{\scriptsize #3}}}
\begin{document}
\title{Nonparametric regression estimation for random fields in a fixed-design}
\author{Mohamed EL MACHKOURI}
\maketitle

{\renewcommand\abstractname{Abstract}
\begin{abstract}
We investigate the nonparametric estimation for regression in a
fixed-design setting when the errors are given by a field of
dependent random variables. Sufficient conditions for kernel
estimators to converge uniformly are obtained. These estimators
can attain the optimal rates of uniform convergence and the
results apply to a large class of random fields which contains
martingale-difference random fields and mixing random fields.
\vspace{8cm}\\
\hspace{-0.7cm}{\em AMS Subject Classifications} (2000): 60G60, 62G08
{\em Short title:} Nonparametric regression in a fixed design.
\end{abstract}
\thispagestyle{empty}
\newpage
\section{Introduction}
Over the last few years nonparametric estimation for random fields
(or spatial processes) was given increasing attention stimulated
by a growing demand from applied research areas (see Guyon
\cite{Guyon}). In fact, spatial data arise in various areas of
research including econometrics, image analysis, meterology, geostatistics...
Our aim in this paper is to investigate uniform strong convergence rates
of a regression estimator in a fixed design
setting when the errors are given by a stationary field of
dependent random variables which show spatial interaction. We
are most interested in conditions which ensure convergence rates
to be identical to those in the case of independent errors (see
Stone \cite{Stone}). Currently the author is working on extensions
of the present results to the random design framework. Let
$\Z^{d}$, $d\geq 1$ denote the integer lattice points in the
$d$-dimensional Euclidean space. By a stationary real random field
we mean any family $(\varepsilon_{k})_{k\in\Z^{d}}$ of real-valued
random variables defined on a probability space $(\Omega, \F, \P)$
such that for any $(k,n)\in\Z^{d}\times\N^{\ast}$ and any
$(i_{1},...,i_{n})\in (\Z^{d})^{n}$, the random vectors
$(\varepsilon_{i_{1}},...,\varepsilon_{i_{n}})$ and
$(\varepsilon_{i_{1}+k},...,\varepsilon_{i_{n}+k})$ have the same
law. The regression model which we are interested in is
\begin{equation}\label{model}
Y_{i}=g(i/n)+\varepsilon_{i},\quad i\in\Lambda_{n}=\{1,...,n\}^{d}
\end{equation}
where $g$ is an unknown smooth function and
$(\varepsilon_{i})_{i\in\Z^{d}}$ is a zero mean stationary real
random field. Note that this model was considered also by Bosq
\cite{Bosq} and Hall et Hart \cite{Hall-Hart} for time series
($d=1$). Let $K$ be a probability kernel defined on $\R^{d}$ and
$(h_{n})_{n\geq 1}$ a sequence of positive numbers which converges
to zero and which satisfies $(nh_{n})_{n\geq 1}$ goes to infinity.
We estimate the function $g$ by the kernel-type estimator $g_{n}$
defined for any $x$ in $[0,1]^{d}$ by
\begin{equation}\label{def-g_n}
g_{n}(x)=\frac{\sum_{i\in\Lambda_{n}}Y_{i}\,
K\left(\frac{x-i/n}{h_{n}}\right)}{\sum_{i\in\Lambda_{n}}
K\left(\frac{x-i/n}{h_{n}}\right)}.
\end{equation}
Note that Assumption $\bold{A1)}$ in section 2 ensures that
$g_{n}$ is well defined. Until now, most of existing theoretical
nonparametric results of dependent random variables pertain to
time series (see Bosq \cite{Bosq-livre}) and relatively few
generalizations to the spatial domain are available. Key
references on this topic are Biau \cite{Biau-2003}, Carbon et al.
\cite{Carbon-Hallin-Tran}, Carbon et al. \cite{Carbon-Tran-Wu},
Hallin et al. \cite{Hallin-Lu-Tran-2001},
\cite{Hallin-Lu-Tran-2004a}, Tran \cite{Tran}, Tran and Yakowitz
\cite{Tran-Yakowitz} and Yao \cite{Yao} who have investigated
nonparametric density estimation for random fields and Altman
\cite{Altman}, Biau and Cadre \cite{Biau-Cadre}, Hallin et al.
\cite{Hallin-Lu-Tran-2004b} and Lu and Chen \cite{Lu-Chen-2002},
\cite{Lu-Chen-2004} who have studied spatial prediction and
spatial regression estimation. The classical asymptotic theory in
statistics is built upon central limit theorems, law of large
numbers and large deviations inequalities for the sequences of
random variables. These classical limit theorems have been
extended to the setting of spatial processes. In particular, some
key results on the central limit theorem and its functional
versions are Alexander and Pyke \cite{Alex-Pyke}, Bass
\cite{Bass}, Basu and Dorea \cite{Bas-Dor}, Bolthausen
\cite{Bolthausen} and more recently Dedecker \cite{JD-tcl},
\cite{JD-tflc}, El Machkouri \cite{MEM-KK} and El Machkouri and
Voln\'y \cite{EM-Volny}. For a survey on limit theorems for
spatial processes and some applications in statistical physics,
one can refer to Nahapetian \cite{Nahapetian-book}. Note also that
the main results (section 3) of this work are obtained via
exponential inequalities for
random fields discovered by El Machkouri \cite{MEM-KK}. \\
The paper is organized as follows. The next section sets up the
notations and the assumptions which will be considered in the
sequel. In section 3, we present our main results on both weak and
strong consistencies rates of the estimator $g_{n}$. The last
section is devoted to the proofs.
\section{Notations and Assumptions}
In the sequel we denote $\|x\|=\max_{1\leq k\leq d}\vert
x_{k}\vert$ for any $x=(x_{1},...,x_{d})\in[0,1]^d$. With a view
to obtain optimal convergence rates for the estimator $g_{n}$
defined by ($\ref{def-g_n}$), we have to make the following
assumptions on the regression function $g$ and the probability
kernel K:
\begin{fleuveA}
\item The probability kernel $K$ is symmetric, nonnegative,
supported by $[-1,1]^{d}$ and satisfies a Lipschitz condition
$\vert K(x)-K(y)\vert\leq\eta\|x-y\|$ for any $x,y\in[-1,1]^{d}$
and some $\eta>0$. In addition there exists $c,C>0$ such that
$c\leq K(x)\leq C$ for any $x\in[-1,1]^{d}$.
\item There exists a constant $B>0$ such that $\vert g(x)-g(y)\vert\leq B\|x-y\|$
for any $x,y\in[0,1]^d$, that is $g$ is $B$-Lipschitz.
\end{fleuveA}
A Young function $\psi$ is a real convex nondecreasing function
defined on $\R^{+}$ which satisfies
$\lim_{t\to\infty}\psi(t)=+\infty$ and $\psi(0)=0$. We define the
Orlicz space $L_{\psi}$ as the space of real random variables $Z$
defined on the probability space $(\Omega, \F, \P)$ such that
$E[\psi(\vert Z\vert/c)]<+\infty$ for some $c>0$. The Orlicz space
$L_{\psi}$ equipped with the so-called Luxemburg norm $\| .
\|_{\psi}$ defined for any real random variable $Z$ by
$$
\| Z\|_{\psi}=\inf\{\,c>0\,;\,E[\psi(\vert Z\vert/c)]\leq 1\,\}
$$
is a Banach space. For more about Young functions and Orlicz
spaces one can refer to Krasnosel'skii and Rutickii \cite{K-R}.
Let $\beta>0$. We denote by $\psi_{\beta}$ the Young function
defined for any $x\in\R^{+}$ by
$$
\psi_{\beta}(x)=\exp((x+\xi_{\beta})^{\beta})-\exp(\xi_{\beta}^{\beta})\quad\textrm{where}\quad
\xi_{\beta}=((1-\beta)/\beta)^{1/\beta}\ind{\{0<\beta<1\}}.
$$
On the lattice $\Z^{d}$ we define the lexicographic order as
follows: if $i=(i_{1},...,i_{d})$ and $j=(j_{1},...,j_{d})$ are
distinct elements of $\Z^{d}$, the notation $i<_{lex}j$ means that
either $i_{1}<j_{1}$ or for some $p$ in $\{2,3,...,d\}$,
$i_{p}<j_{p}$ and $i_{q}=j_{q}$ for $1\leq q<p$. Let the sets
$\{V_{i}^{k}\,;\,i\in\Z^{d}\,,\,k\in\N^{\ast}\}$ be defined as
follows:
$$
V_{i}^{1}=\{j\in\Z^{d}\,;\,j<_{lex}i\},
$$
and for $k\geq 2$
$$
V_{i}^{k}=V_{i}^{1}\cap\{j\in\Z^{d}\,;\,\vert i-j\vert\geq
k\}\quad\textrm{where}\quad \vert i-j\vert=\max_{1\leq l\leq
d}\vert i_{l}-j_{l}\vert.
$$
For any subset $\Gamma$ of $\Z^{d}$ define
$\F_{\Gamma}=\sigma(\varepsilon_{i}\,;\,i\in\Gamma)$ and set
$$
E_{\vert
k\vert}(\varepsilon_{i})=E(\varepsilon_{i}\vert\F_{V_{i}^{\vert
k\vert}}),\quad k\in V_{i}^{1}.
$$
Denote $\beta(q)=2q/(2-q)$ for $0<q<2$ and consider the following
conditions:
\begin{fleuveC}
\item $\varepsilon_{0}\in L^{\infty}$ and
$$
\sum_{k\in V_{0}^{1}}\|\varepsilon_{k}E_{\vert
k\vert}(\varepsilon_{0})\|_{\infty}<\infty.
$$
\item There exists $0<q<2$ such that $\varepsilon_{0}\in L_{\psi_{\beta(q)}}$ and
$$
\sum_{k\in V_{0}^{1}}\left\|\sqrt{\vert \varepsilon_{k}E_{\vert
k\vert}(\varepsilon_{0})\vert}\right\|^{2}_{\psi_{\beta(q)}}<\infty.
$$
\item There exists $p>2$ such that $\varepsilon_{0}\in L^{p}$ and
$$
\sum_{k\in V_{0}^{1}}\|\varepsilon_{k}E_{\vert
k\vert}(\varepsilon_{0})\|_{\frac{p}{2}}<\infty.
$$
\item $\varepsilon_{0}\in L^{2}$ and
$\sum_{k\in\Z^{d}}\vert
E(\varepsilon_{0}\varepsilon_{k})\vert<\infty$.
\end{fleuveC}
\begin{RMQ}
{\em Note that Dedecker \cite{JD-tcl} established the central
limit theorem for any stationary square-integrable random field
$(\varepsilon_{k})_{k\in\Z^{d}}$ which satisfies the condition
$\sum_{k\in V_{0}^{1}}\|\varepsilon_{k}E_{\vert
k\vert}(\varepsilon_{0})\|_{1}<\infty$.}
\end{RMQ}
In classical statistical physics, there exists spatial processes
which satisfy conditions
\emph{$\bold{C1)}$},...,\emph{$\bold{C4)}$}. For example,
Nahapetian and Petrosian \cite{Nahapetian-Petrosian} gave
sufficient conditions for a Gibbs field
$(\varepsilon_{k})_{k\in\Z^{d}}$ to possess the following
martingale difference property: for any $i$ in $\Z^d$,
$E(\varepsilon_{i}\vert\F_{V_{i}^1})=0$ a.s. Another examples of
random fields which satisfy conditions
\emph{$\bold{C1)}$},...,\emph{$\bold{C4)}$} can be found also
among the class of mixing random fields. More precisely, given two
sub-$\sigma$-algebras $\U$ and $\V$ of $\F$, different measures of
their dependence have been considered in the literature. We are
interested by two of them. The $\alpha$-mixing and $\phi$-mixing
coefficients had been introduced by Rosenblatt \cite{Ros} and
Ibragimov \cite{Ibrag} respectively and can be defined by
\begin{align*}
\alpha(\U,\V)&=\sup\{\vert\P(U\cap V)-\P(U)\P(V)\vert,\,U\in\U
,\,V\in\V\}\\
\phi(\U,\V)&=\sup\{\| \P(V\vert \U)-\P(V)\|_{\infty}\,
,\,V\in\V\}.
\end{align*}
We have $2\alpha(\U,\V)\leq\phi(\U,\V)$ and these coefficients
equal zero if and only if the $\sigma$-algebras $\U$ and $\V$ are
independent. Denote by $\sharp\Gamma$ the cardinality of any
subset $\Gamma$ of $\Z^{d}$. In the sequel, we shall use the
following non-uniform mixing coefficients defined for any
$(k,l,n)$ in $(\N^{\ast}\cup\{\infty\})^{2}\times\N$ by
\begin{align*}
\alpha_{k,l}(n)&=\sup\,\{\alpha(\F_{\Gamma_{1}},\F_{\Gamma_{2}}),\,
\sharp\Gamma_{1}\leq k,\, \sharp\Gamma_{2}\leq l,\,
\rho(\Gamma_{1},\Gamma_{2})\geq n\},\\
\phi_{k,l}(n)&=\sup\,\{\phi(\F_{\Gamma_{1}},\F_{\Gamma_{2}}),\,
\sharp\Gamma_{1}\leq k,\, \sharp\Gamma_{2}\leq l,\,
\rho(\Gamma_{1},\Gamma_{2})\geq n\},
\end{align*}
where the distance $\rho$ is defined by
$\rho(\Gamma_{1},\Gamma_{2})=\min\{\vert
i-j\vert,\,i\in\Gamma_{1},\,j\in\Gamma_{2}\}$. We say that the
random field $(\varepsilon_{k})_{k\in\Z^{d}}$ is $\alpha$-mixing
or $\phi$-mixing if there exists a pair $(k,l)$ in
$(\N^{\ast}\cup\{\infty\})^{2}$ such that $\lim_{n\to
\infty}\alpha_{k,l}(n)=0$ or $\lim_{n\to\infty}\phi_{k,l}(n)=0$
respectively. For more about mixing coefficients one can refer to
Doukhan \cite{Doukhan}. We consider the following mixing conditions:\\
\\
\textbf{$\bold{C^{'}1)}$} $\varepsilon_{0}\in
L^{\infty}$ and
$$
\sum_{k\in\Z^{d}}\phi_{\infty,1}(\vert k\vert)<\infty.
$$
\textbf{$\bold{C^{'}2)}$} There exists $0<q<2$ such that $\varepsilon_{0}\in
L_{\psi_{\beta(q)}}$ and
$$
\sum_{k\in\Z^{d}}\sqrt{\phi_{\infty,1}(\vert k\vert)}<\infty
$$
or
$$
\sum_{k\in\Z^{d}}c_{k}^{2}(\beta(q))<\infty
$$
where for any $\beta>0$
\begin{equation}\label{definition-ck}
c_{k}(\beta)=\inf\left\{c>0\,\big\vert\,\int_{0}^{\alpha_{1,\infty}(\vert
k\vert)}\psi_{\beta}\left(\frac{Q_{\varepsilon_{0}}(u)}{c}\right)\,du\leq
1\right\}.
\end{equation}
\textbf{$\bold{C^{'}3)}$} There exists $p>2$ such that
$\varepsilon_{0}\in L^{p}$ and
\begin{equation}\label{mixing-quantile}
\sum_{k\in\Z^{d}}\left(\int_{0}^{\alpha_{1,\infty}(\vert
k\vert)}Q_{\varepsilon_{0}}^{p}(u)\,du\right)^{2/p}<\infty
\end{equation}
where $Q_{\varepsilon_{0}}$ is the inverse cadlag of the tail
function $t\to\P(\vert \varepsilon_{0}\vert>t)$ (i.e. for any
$u\geq 0$,
$Q_{\varepsilon_{0}}(u)=\inf\left\{t>0\,\vert\,\P(\vert\varepsilon_{0}\vert>t)\leq
u\right\})$.
\begin{RMQ}{\em Let us note that if $p=2+\delta$ for some $\delta>0$ then the
condition
$$
\sum_{m=1}^{\infty}m^{d-1}\alpha_{1,\infty}^{\frac{\delta}{2+\delta}-\varepsilon}(m)<\infty\quad\textrm{for
some $\varepsilon>0$}
$$
is more restrictive than condition ($\ref{mixing-quantile}$) and
is known to be sufficient for the random field
$(\varepsilon_{k})_{k\in\Z^{d}}$ to satisfy a functional central
limit theorem (cf. Dedecker \cite{JD-tflc}).}
\end{RMQ}
In statistical physics, using the Dobrushin's uniqueness condition
(cf. \cite{Dobrushin}), one can construct Gibbs fields satisfying
a uniform exponential mixing condition which is more restrictive
than conditions \emph{$\bold{C^{'}1)}$}, \emph{$\bold{C^{'}2)}$}
and \emph{$\bold{C^{'}3)}$} (see Guyon \cite{Guyon}, theorem
2.1.3, p. 52).
\section{Main results}
Let $(Z_{n})_{n\geq 1}$ be a sequence of real random variables and
$(v_{n})_{n\geq 1}$ be a sequence of positive numbers. We say that
$$
Z_{n}=O_{a.s.}\left[v_{n}\right]
$$
if there exists $\lambda>0$ such that
$$
\limsup_{n\to\infty}\frac{\vert
Z_{n}\vert}{v_{n}}\leq\lambda\quad\textrm{a.s.}
$$
Our main result is the following.
\begin{Th}\label{rate-variance-ps} Assume that the assumption
$\emph{\textbf{A1)}}$ holds.
\begin{fleuve}
\item If \emph{$\bold{C1)}$} holds then
\begin{equation}\label{rvps3}
\sup_{x\in [0,1]^{d}}\vert
g_{n}(x)-Eg_{n}(x)\vert=O_{a.s.}\left[\frac{(\log
n)^{1/2}}{(nh_{n})^{d/2}}\right].
\end{equation}
\item If \emph{$\bold{C2)}$} holds for some $0<q<2$ then
\begin{equation}\label{rvps2}
\sup_{x\in [0,1]^{d}}\vert
g_{n}(x)-Eg_{n}(x)\vert=O_{a.s.}\left[\frac{(\log
n)^{1/q}}{(nh_{n})^{d/2}}\right].
\end{equation}
\item Assume that \emph{$\bold{C3)}$} holds for some $p>2$ and
$h_{n}=n^{-\theta_{2}}(\log n)^{\theta_{1}}$ for some
$\theta_{1},\theta_{2}>0$. Let $a,b\geq 0$ be fixed and denote
$$
v_{n}=\frac{n^{a}(\log
n)^{b}}{(nh_{n})^{d/2}}\quad\textrm{and}\quad\theta=\frac{2a(d+p)-d^{2}-2}{d(3d+2)}.
$$
If $\theta\geq\theta_{2}$ and $d(3d+2)\theta_{1}+2(d+p)b>2$ then
\begin{equation}\label{rvps1}
\sup_{x\in[0,1]^{d}}\vert
g_{n}(x)-Eg_{n}(x)\vert=O_{a.s.}\left[v_{n}\right].
\end{equation}
\end{fleuve}
\end{Th}
\begin{RMQ}{\em Theorem $\ref{rate-variance-ps}$ shows that the optimal uniform convergence
rate is obtained for bounded errors (cf. estimation ($\ref{rvps3}$)) and that it is ``almost'' optimal if one considers errors with only finite
exponential moments (cf. estimation ($\ref{rvps2}$)).}
\end{RMQ}
\begin{Th}\label{rate-variance-proba} Assume that the assumption
$\emph{\textbf{A1)}}$ holds.
\begin{fleuve}
\item Assume that \emph{$\bold{C3)}$} holds for some $p>2$. Let
$a>0$ be fixed and denote
$$
v_{n}=\frac{n^{a}}{(nh_{n})^{d/2}}\quad\textrm{and}\quad\theta=\frac{2a(d+p)-d^{2}}{d(3d+2)}.
$$
If $\theta>0$ and $h_{n}\geq n^{-\theta}$ then
\begin{equation}\label{rvLp}
\left\|\sup_{x\in[0,1]^{d}}\vert
g_{n}(x)-Eg_{n}(x)\vert\right\|_{p}=O\left[v_{n}\right].
\end{equation}
\item If \emph{$\bold{C4)}$} holds then
\begin{equation}\label{rvL2}
\sup_{x\in[0,1]^{d}}\left\|g_{n}(x)-Eg_{n}(x)\right\|_{2}
=O\left[(nh_{n})^{-d/2}\right].
\end{equation}
\end{fleuve}
\end{Th}
In the sequel, we denote by $\textrm{Lip}(B)$ the set of
$B$-Lipschitz functions. The following proposition gives the
convergence of $Eg_{n}(x)$ to $g(x)$.
\begin{Prop}\label{rate-biais}
Assume that the assumption $\emph{\textbf{A2)}}$ holds then
$$
\sup_{x\in [0,1]^{d}}\sup_{g\in\,\textrm{Lip}(B)}\vert
Eg_{n}(x)-g(x)\vert=O\left[h_{n}\right].
$$
\end{Prop}
From Proposition $\ref{rate-biais}$ and Theorem
$\ref{rate-variance-ps}$ we derive the following corollary.
\begin{Cor}\label{cor-optimal-rate-ps}
Assume that $\emph{\textbf{A1)}}$ and $\emph{\textbf{A2)}}$ hold
and let $h_{n}=\left(n^{-d}\log n\right)^{1/(2+d)}$.
\begin{fleuve}
\item If \emph{$\bold{C1)}$} holds then
\begin{equation}\label{optimal-rps3}
\sup_{x\in [0,1]^{d}}\sup_{g\in\,\textrm{Lip}(B)}\vert
g_{n}(x)-g(x)\vert=O_{a.s.}\left[\left(\frac{\log
n}{n^d}\right)^{\frac{1}{2+d}}\right].
\end{equation}
\item If \emph{$\bold{C2)}$} holds for some $0<q<2$ then
\begin{equation}\label{optimal-rps2}
\sup_{x\in [0,1]^{d}}\sup_{g\in\,\textrm{Lip}(B)}\vert
g_{n}(x)-g(x)\vert=O_{a.s.}\left[u(n)\left(\frac{\log
n}{n^d}\right)^{\frac{1}{2+d}}\right]
\end{equation}
where $u(n)=(\log n)^{(2-q)/2q}$.
\item Let $\varepsilon>0$ be fixed. If \emph{$\bold{C3)}$} holds for some $p>2$ satisfying
\begin{equation}\label{min-p}
p\geq\frac{4d^3+(4-2\varepsilon)d^2+(2-4\varepsilon)d+4}{2\varepsilon(2+d)}
\end{equation}
then
\begin{equation}\label{optimal-rps1}
\sup_{x\in [0,1]^{d}}\sup_{g\in\,\textrm{Lip}(B)}\vert
g_{n}(x)-g(x)\vert=O_{a.s.}\left[u(n)\left(\frac{\log
n}{n^d}\right)^{\frac{1}{2+d}}\right]
\end{equation}
where $u(n)=n^{\varepsilon}$.
\end{fleuve}
\end{Cor}
\begin{RMQ}{\em Note that the consistency rate $(n^{-d}\log n)^{1/(2+d)}$
is known to be the optimal one (see Stone \cite{Stone}).}
\end{RMQ}
From Proposition $\ref{rate-biais}$ and Theorem
$\ref{rate-variance-proba}$ we derive the following corollary.
\begin{Cor}\label{cor-optimal-rate-proba}
Assume that $\emph{\textbf{A1)}}$ and $\emph{\textbf{A2)}}$ hold
and let $h_{n}=n^{-d/(2+d)}$.
\begin{fleuve}
\item Let $\varepsilon>0$ be fixed. If \emph{$\bold{C3)}$} holds for some $p>2$ satisfying
\begin{equation}\label{min-p-bis}
p\geq\frac{4d^3+(4-2\varepsilon)d^2-4\varepsilon
d}{2\varepsilon(2+d)}
\end{equation}
then
\begin{equation}\label{optimal-rproba2}
\left\|\sup_{x\in [0,1]^{d}}\sup_{g\in\,\textrm{Lip}(B)}\vert
g_{n}(x)-g(x)\vert\right\|_{p}=O\left[n^{-\frac{d}{2+d}+\varepsilon}\right].
\end{equation}
\item If \emph{$\bold{C4)}$} holds then
\begin{equation}\label{optimal-rproba1}
\sup_{x\in[0,1]^{d}}\left\|\sup_{g\in\,\textrm{Lip}(B)}\vert
g_{n}(x)-g(x)\vert\right\|_{2}=O\left[n^{-\frac{d}{2+d}}\right].
\end{equation}
\end{fleuve}
\end{Cor}
Finally the rates of convergence obtained above are valid when the
errors are given by a mixing random field. More precisely, we have
the following corollary.
\begin{Cor}\label{cor-optimal-rate-ps-mixing}
Theorems $\ref{rate-variance-ps}$ and $\ref{rate-variance-proba}$
and Corollaries $\ref{cor-optimal-rate-ps}$ and
$\ref{cor-optimal-rate-proba}$ still hold if one replace
conditions \emph{$\bold{C1)}$}, \emph{$\bold{C2)}$} and
\emph{$\bold{C3)}$} by conditions \emph{$\bold{C^{'}1)}$},
\emph{$\bold{C^{'}2)}$} and \emph{$\bold{C^{'}3)}$} respectively.
\end{Cor}
\section{Proofs}
For any $x$ in $[0,1]^{d}$ and any integer $n\geq 1$ we define
$B_{n}(x)=Eg_{n}(x)-g(x)$ and $V_{n}(x)=g_{n}(x)-Eg_{n}(x)$. More
precisely
\begin{align*}
B_{n}(x)&=\frac{\sum_{i\in\Lambda_{n}}a_{i}(x)g(i/n)}{\sum_{i\in\Lambda_{n}}a_{i}(x)}-g(x)\\
V_{n}(x)&=\frac{\sum_{i\in\Lambda_{n}}a_{i}(x)\varepsilon_{i}}{\sum_{i\in\Lambda_{n}}a_{i}(x)}
\end{align*}
where $a_{i}(x)=K\left(\frac{x-i/n}{h_{n}}\right)$. In the sequel,
we denote also
$S_{n}(x)=\sum_{i\in\Lambda_{n}}a_{i}(x)\varepsilon_{i}$ for any
$x\in[0,1]^{d}$. We start with the following lemma.
\begin{lemma}\label{lemma-borne-de-K}
There exists constants $c,C>0$ such that for any $x\in [0,1]^{d}$
and any $n\in\N^{\ast}$,
\begin{equation}\label{borne-de-K}
c\prod_{k=1}^{d}[n(x_{k}+h_{n})]\leq\sum_{i\in\Lambda_{n}}a_{i}(x)\leq
C\prod_{k=1}^{d}[n(x_{k}+h_{n})]
\end{equation}
where $[\,.\,]$ denote the integer part function.
\end{lemma}
{\em Proof of Lemma $\ref{lemma-borne-de-K}$.} Since the kernel
$K$ is supported by $[-1,1]^d$, we have
$$
\sum_{i\in\Lambda_{n}}a_{i}(x)
=\sum_{i_{1}=1}^{[n(x_{1}+h_{n})]}\ldots\sum_{i_{d}=1}^{[n(x_{d}+h_{n})]}a_{i}(x).
$$
By assumption,
there exists constants $c,C>0$ such that $c\leq K(y)\leq C$ for
any
$y\in [-1,1]^{d}$. The proof of Lemma $\ref{lemma-borne-de-K}$ is complete.
\subsection{Proof of Theorem $\textbf{\ref{rate-variance-ps}}$}
Let $(\upsilon_{n})_{n\geq 1}$ be a sequence of positive numbers
going to zero. Following Carbon and al. \cite{Carbon-Tran-Wu} the
compact set $[0,1]^{d}$ can be covered by $r_{n}$ cubes $I_{k}$
having sides of length $l_{n}=\upsilon_{n}h_{n}^{2d+1}$ and center
at $c_{k}$. Clearly there exists $c>0$ such that $r_{n}\leq
c/l_{n}^{d}$. Define
\begin{align*}
A_{1,n}(g)&=\max_{1\leq k\leq r_{n}}\sup_{x\in I_{k}}\vert
g_{n}(x)-g_{n}(c_{k})\vert\\
A_{2,n}(g)&=\max_{1\leq k\leq r_{n}}\sup_{x\in I_{k}}\vert
Eg_{n}(x)-Eg_{n}(c_{k})\vert\\
A_{3,n}&=\max_{1\leq k\leq r_{n}}\vert
g_{n}(c_{k})-Eg_{n}(c_{k})\vert
\end{align*}
then
\begin{equation}\label{inegalite}
\sup_{x\in [0,1]^{d}}\vert
g_{n}(x)-Eg_{n}(x)\vert\leq\sup_{g\in\,\textrm{Lip}(B)}[A_{1,n}(g)+A_{2,n}(g)]+A_{3,n}.
\end{equation}
\begin{lemma}\label{A1n-A2n} For $i=1,2$ we have
$$
\sup_{g\in\,\textrm{Lip}(B)}A_{i,n}(g)=O_{a.s.}\left[v_{n}\right].
$$
\end{lemma}
{\em Proof of Lemma $\ref{A1n-A2n}$.} Since
$g\in\,\textrm{Lip}(B)$, we can assume without loss of generality
that $g$ is bounded by $B$ on the set $[0,1]^d$. For any $x\in
I_{k}$, we have
$$
g_{n}(x)-g_{n}(c_{k})=\sigma_{1}+\sigma_{2}
$$
where
$$
\sigma_{1}=\frac{\sum_{i\in\Lambda_{n}}Y_{i}(a_{i}(x)-a_{i}(c_{k}))}{\sum_{i\in\Lambda_{n}}a_{i}(x)}
$$
and
$$
\sigma_{2}=\frac{\sum_{i\in\Lambda_{n}}(a_{i}(c_{k})-a_{i}(x))}
{\sum_{i\in\Lambda_{n}}a_{i}(x)\times\sum_{i\in\Lambda_{n}}a_{i}(c_{k})}
\sum_{i\in\Lambda_{n}}Y_{i}a_{i}(c_{k}).
$$
Now, by Lemma $\ref{lemma-borne-de-K}$ and Assumption
$\bold{A1)}$, we derive that there exists constants $c,\eta>0$
such that for any $n$ sufficiently large
$$
\vert\sigma_{1}\vert\leq \frac{2^d\eta
l_{n}/h_{n}}{c(nh_{n})^{d}}\sum_{i\in\Lambda_{n}}\vert
Y_{i}\vert\leq\frac{\eta
v_{n}h_{n}^{d}}{c}\left(B+\frac{1}{n^{d}}\sum_{i\in\Lambda_{n}}\vert\varepsilon_{i}\vert\right)
$$
and
$$
\vert\sigma_{2}\vert\leq\frac{4^d\eta
n^{d}l_{n}/h_{n}}{c^{2}(nh_{n})^{2d}}\sum_{i\in\Lambda_{n}}\vert
Y_{i}\vert\leq\frac{\eta
v_{n}}{c^{2}}\left(B+\frac{1}{n^{d}}\sum_{i\in\Lambda_{n}}\vert\varepsilon_{i}\vert\right)
$$
Since $(\varepsilon_{i})$ is a stationary ergodic random field the
lemma easily follows from the last inequalities and the Birkhoff
ergodic theorem. The proof of Lemma $\ref{A1n-A2n}$ is complete.
\begin{lemma}\label{A3n}
Assume that either \emph{$\bold{C1)}$} holds and $v_{n}=(\log
n)^{1/2}/(nh_{n})^{d/2}$ or \emph{$\bold{C2)}$} holds for some
$0<q<2$ and $v_{n}=(\log n)^{1/q}/(nh_{n})^{d/2}$ then
$$
A_{3,n}=O_{a.s.}[v_{n}]
$$
\end{lemma}
\emph{Proof of Lemma $\ref{A3n}$.} Let $0<q\leq 2$ be fixed. We
consider the exponential Young function define for any
$x\in\R^{+}$ by
$\psi_{q}(x)=\exp((x+\xi_{q})^{q})-\exp(\xi_{q}^{q})$ where
$\xi_{q}=((1-q)/q)^{1/q}\ind{\{0<q<1\}}$. Let $\lambda>0$ and
$x\in [0,1]^{d}$ be fixed
\begin{align*}
\P\left(\vert V_{n}(x)\vert>\lambda v_{n}\right)
&=\P\left(\big\vert S_{n}(x)\big\vert>\lambda v_{n}\sum_{i\in\Lambda_{n}}a_{i}(x)\right)\\
&\leq (1+e^{\xi_{q}^{q}})\,\exp\bigg[-\left(\frac{\lambda
v_{n}\sum_{i\in\Lambda_{n}}a_{i}(x)}{\vert\vert\sum_{i\in\Lambda_{n}}a_{i}(x)\varepsilon_{i}\vert\vert_{\psi_{q}}}+\xi_{q}\right)^{q}\bigg].
\end{align*}
For any $i\in\Lambda_{n}$ and any $0<q<2$ denote
\begin{equation}\label{biq}
b_{i,q}(a(x)\varepsilon)=
\big\|a_{i}(x)\varepsilon_{i}\big\|_{\psi_{\beta(q)}}^{2}+\sum_{k\in
V_{i}^{1}}\left\|\sqrt{\big\vert a_{k}(x)\varepsilon_{k}E_{\vert
k-i\vert}(a_{i}(x)\varepsilon_{i})\big\vert}\right\|_{\psi_{\beta(q)}}^{2}
\end{equation}
and
\begin{equation}\label{bi2}
b_{i,2}(a(x)\varepsilon)=
\big\|a_{i}(x)\varepsilon_{i}\big\|_{\infty}^{2}+\sum_{k\in
V_{i}^{1}}\| a_{k}(x)\varepsilon_{k}E_{\vert
k-i\vert}(a_{i}(x)\varepsilon_{i})\|_{\infty}
\end{equation}
where $V_{i}^{1}=\{j\in\Z^{d}\,;\,j<_{lex}i\}$. Using
Kahane-Khintchine inequalities (cf. El Machkouri \cite{MEM-KK},
Theorem 1) we derive that if Condition \emph{$\bold{C2)}$} holds
for some $0<q<2$ then
\begin{equation}\label{inegalite5}
\P\left(\vert V_{n}(x)\vert>\lambda v_{n}\right)\leq
(1+e^{\xi_{q}^{q}})\,\exp\bigg[-\left(\frac{\lambda\,v_{n}\sum_{i\in\Lambda_{n}}a_{i}(x)}
{M(\sum_{i\in\Lambda_{n}}b_{i,q}(a(x)\varepsilon))^{1/2}}+\xi_{q}\right)^{q}\bigg]
\end{equation}
where $M$ is a positive constant depending only on $q$ and on the
probability kernel $K$. Now using the definition ($\ref{biq}$) and
Lemma $\ref{lemma-borne-de-K}$ there exist constants $c,M>0$ such
that
\begin{align*}
\sup_{x\in[0,1]^{d}}\P\left(\vert V_{n}(x)\vert>\lambda
v_{n}\right)&\leq
(1+e^{\xi_{q}^{q}})\,\exp\bigg[-\left(\frac{\lambda\,v_{n}\,(\sum_{i\in\Lambda_{n}}a_{i}(x))^{1/2}}
{M}+\xi_{q}\right)^{q}\bigg]\\
&\leq(1+e^{\xi_{q}^{q}})\,\exp\bigg[-\frac{c^q\,\lambda^{q}\,v_{n}^{q}\,([nh_{n}])^{dq/2}}
{M^{q}}\bigg]
\end{align*}
So if $v_{n}=(\log n)^{1/q}/(nh_{n})^{d/2}$ and $n$ is sufficiently large then
\begin{equation}\label{inegalite6}
\sup_{x\in[0,1]^{d}}\P\left(\vert V_{n}(x)\vert>\lambda
v_{n}\right)\leq(1+e^{\xi_{q}^{q}})\,\exp\bigg[-\frac{c^q\,\lambda^{q}\,\log
n} {2^{dq/2}M^{q}}\bigg].
\end{equation}
If Condition \emph{$\bold{C1)}$} holds then ($\ref{inegalite5}$)
still hold with $q=2$ (cf. El Machkouri \cite{MEM-KK}, Theorem 1).
So if $v_{n}=(\log n)^{1/2}/(nh_{n})^{d/2}$ and $n$ is large it follows that
\begin{equation}\label{inegalite7}
\sup_{x\in[0,1]^{d}}\P\left(\vert V_{n}(x)\vert>\lambda
v_{n}\right)\leq 2\exp\bigg[-\frac{c^2\,\lambda^{2}\,\log
n}{2^{d}M^{2}}\bigg].
\end{equation}
Since
$$
\P\left(\vert A_{3,n}\vert>\lambda v_{n}\right) \leq
r_{n}\sup_{x\in[0,1]^{d}}\P\left(\vert V_{n}(x)\vert>\lambda
v_{n}\right),
$$
using ($\ref{inegalite6}$) and ($\ref{inegalite7}$), choosing
$\lambda$ sufficiently large and applying Borel-Cantelli's lemma,
we derive
$$
\P\left(\limsup_{n\to\infty}\{\vert A_{3,n}\vert>\lambda
v_{n}\}\right)=0
$$
and
$$
\P\left(\limsup_{n\to\infty}\frac{\vert
A_{3,n}\vert}{v_{n}}\leq\lambda\right)=1.
$$
The proof of points $\bold{1)}$ and $\bold{2)}$ of Theorem
$\ref{rate-variance-ps}$ are completed by combining Inequality
($\ref{inegalite}$) with Lemmas $\ref{A1n-A2n}$ and $\ref{A3n}$.
\begin{lemma}\label{A3n-nonborne}
Assume that \emph{$\bold{C3)}$} holds for some $p>2$ and
$h_{n}=n^{-\theta_{2}}(\log n)^{\theta_{1}}$ for some
$\theta_{1},\theta_{2}>0$. Let $a,b\geq 0$ be fixed and denote
$$
v_{n}=\frac{n^{a}(\log
n)^{b}}{(nh_{n})^{d/2}}\quad\textrm{and}\quad\theta=\frac{2a(d+p)-d^{2}-2}{d(3d+2)}.
$$
If $\theta\geq\theta_{2}$ and $d(3d+2)\theta_{1}+2(d+p)b>2$ then
$$
\lim_{n\to+\infty}\frac{\vert
A_{3,n}\vert}{v_{n}}=0\quad\textrm{a.s.}
$$
\end{lemma}
{\em Proof of Lemma $\ref{A3n-nonborne}$.} Let $p>2$ be fixed. For
any $\lambda>0$
\begin{align*}
\P\left(\vert V_{n}(x)\vert>\lambda v_{n}\right)
&=\P\left(\vert S_{n}(x)\vert>\lambda v_{n}\sum_{i\in\Lambda_{n}}a_{i}(x)\right)\\
&\leq\frac{\lambda^{-p}E\vert S_{n}(x)\vert^{p}}{v_{n}^{p}(\sum_{i\in\Lambda_{n}}a_{i}(x))^{p}}\\
&\leq\frac{\lambda^{-p}}{v_{n}^{p}(\sum_{i\in\Lambda_{n}}a_{i}(x))^{p}}\left(2p\sum_{i\in\Lambda_{n}}c_{i}(x)\right)^{p/2}
\end{align*}
where
$c_{i}(x)=a_{i}(x)^{2}\|\varepsilon_{i}\|_{p}^{2}+a_{i}(x)\sum_{k\in
V_{i}^{1}}a_{k}(x)\|\varepsilon_{k}E_{\vert
k-i\vert}(\varepsilon_{i})\|_{\frac{p}{2}}$. The last estimate
follows from a Marcinkiewicz-Zygmund type inequality by Dedecker
(see \cite{JD-tflc}) for real random fields. Noting that there
exists $\gamma>0$ such that $c_{i}(x)\leq\gamma
a_{i}(x),\,x\in[0,1]^{d}$ and using Lemma
$\ref{lemma-borne-de-K}$, we derive that there exists
$\gamma^{'}>0$ such that
$$
\P\left(\vert A_{3,n}\vert>\lambda v_{n}\right)
\leq r_{n}\sup_{x\in[0,1]^{d}}\P\left(\vert V_{n}(x)\vert>\lambda
v_{n}\right)\leq \frac{\gamma^{'}}{\tau_{n}\lambda^{p}}
$$
where $\tau_{n}=l_{n}^{d}v_{n}^{p}([nh_{n}])^{dp/2}$. Since
$v_{n}=n^{a}(\log n)^{b}/(nh_{n})^{d/2}$ and
$l_{n}=v_{n}h_{n}^{2d+1}$ it follows
$$
\frac{1}{\tau_{n}}=\frac{(nh_{n})^{d(d+p)/2}}{h_{n}^{d(2d+1)}n^{a(d+p)}(\log
n)^{b(d+p)}([nh_{n}])^{dp/2}}.
$$
If $n$ is sufficiently large, we derive
\begin{align*}
\frac{1}{\tau_{n}}&\leq\frac{2^{dp/2}(nh_{n})^{d(d+p)/2}}{h_{n}^{d(2d+1)}n^{a(d+p)}(\log n)^{b(d+p)}(nh_{n})^{dp/2}}\\
&=\frac{2^{dp/2}}{h_{n}^{d(3d+2)/2}n^{a(d+p)-d^{2}/2}(\log n)^{b(d+p)}}\\
&\leq\frac{2^{dp/2}}{n(\log
n)^{b(d+p)+\theta_{1}d(3d+2)/2}}\quad\textrm{since $\theta\geq
\theta_{2}$}.
\end{align*}
Now $b(d+p)+\theta_{1}d(3d+2)/2>1$ implies $\sum_{n\geq
1}\tau_{n}^{-1}<\infty$. Applying Borel-Cantelli's lemma, it
follows that for any $\lambda>0$
$$
\P\left(\limsup_{n\to\infty}\{\vert A_{3,n}\vert>\lambda
v_{n}\}\right)=0,
$$
that is for any $\lambda>0$
$$
\P\left(\limsup_{n\to\infty}\frac{\vert
A_{3,n}\vert}{v_{n}}\leq\lambda\right)=1.
$$
The proof of Lemma $\ref{A3n-nonborne}$ is complete and the point
$\bold{3)}$ of Theorem $\ref{rate-variance-ps}$ is obtained by
combining Inequality ($\ref{inegalite}$) with Lemmas
$\ref{A1n-A2n}$ and $\ref{A3n-nonborne}$. The proof of Theorem
$\ref{rate-variance-ps}$ is complete.
\subsection{Proof of Theorem $\textbf{\ref{rate-variance-proba}}$}
We follow the first part of the proof of Theorem
$\ref{rate-variance-ps}$ and we consider the estimation
($\ref{inegalite}$).
\begin{lemma}\label{A3n-nonborne-proba}
Assume that \emph{$\bold{C3)}$} holds for some $p>2$. Let $a>0$ be
fixed and denote
$$
v_{n}=\frac{n^{a}}{(nh_{n})^{d/2}}\quad\textrm{and}\quad\theta=\frac{2a(d+p)-d^{2}}{d(3d+2)}.
$$
If $\theta>0$ and $h_{n}\geq n^{-\theta}$ then
$$
\left\|A_{3,n}\right\|_{p}=O\left[v_{n}\right].
$$
\end{lemma}
{\em Proof of Lemma $\ref{A3n-nonborne-proba}$.} Let $p>2$ and
$x\in[0,1]^d$ be fixed. Using the Marcinkiewicz-Zygmund type
inequality by Dedecker (see \cite{JD-tflc}) as in the proof of
Lemma $\ref{A3n-nonborne}$ there exist $\gamma^{''},c>0$ such that
\begin{align*}
\|V_{n}(x)\|_{p}
&=\left(\frac{E\vert S_{n}(x)\vert^p}{\left(\sum_{i\in\Lambda_{n}}a_{i}(x)\right)^p}\right)^{1/p}\\
&\leq\gamma^{''}\left(\sum_{i\in\Lambda_{n}}a_{i}(x)\right)^{-1/2}\\
&\leq\frac{\gamma^{''}}{\sqrt{c}}([nh_{n}])^{-d/2}\qquad\textrm{by
Lemma $\ref{lemma-borne-de-K}$}.
\end{align*}
It follows that
$$
r_{n}^{1/p}\sup_{x\in[0,1]^{d}}\|V_{n}(x)\|_{p}=O\left[\frac{v_{n}}{\tau_{n}}\right]
$$
where $\tau_{n}=l_{n}^{d/p}v_{n}([nh_{n}])^{d/2}$. If $n$ is sufficiently large then $\tau_{n}\geq 2^{-d/2}l_{n}^{d/p}v_{n}(nh_{n})^{d/2}$, hence using $h_{n}\geq
n^{-\theta}$ we obtain $\tau_{n}\geq 2^{-d/2}$. Finally, we derive
$$
\|A_{3,n}\|_{p}=\|\max_{1\leq k\leq r_{n}}\vert
V_{n}(x_{k})\vert\|_{p}\leq
r_{n}^{1/p}\sup_{x\in[0,1]^{d}}\|V_{n}(x)\|_{p}=O\left[v_{n}\right].
$$
The proof of Lemma $\ref{A3n-nonborne-proba}$ is complete. The
point $\bold{1)}$ of Theorem $\ref{rate-variance-proba}$ is
obtained by combining inequality
($\ref{inegalite}$) and lemmas $\ref{A1n-A2n}$ and $\ref{A3n-nonborne-proba}$.\\
\\
Now, we are going to prove the point $\bold{2)}$ of Theorem
$\ref{rate-variance-proba}$. We have
\begin{align*}
E(S_{n}(x)^{2})&=\sum_{k,l\in\Lambda_{n}}a_{k}(x)\,a_{l}(x)E(\varepsilon_{k}\varepsilon_{l})\\
&=\sum_{k\in\Lambda_{n}}a_{k}(x)^{2}E(\varepsilon_{k}^{2})+\sum_{k\neq
l}a_{k}(x)a_{l}(x)E(\varepsilon_{k}\varepsilon_{l})\\
&=E(\varepsilon_{0}^{2})\sum_{k\in\Lambda_{n}}a_{k}(x)^{2}
+\sum_{k\in\Lambda_{n}}a_{k}(x)\sum_{l\in\Lambda_{n}\backslash\{k\}}a_{l}(x)E(\varepsilon_{k}\varepsilon_{l})\\
&\leq\sum_{l\in\Z^{d}}\vert
E(\varepsilon_{0}\varepsilon_{l})\vert\times\sum_{k\in\Lambda_{n}}a_{k}(x).
\end{align*}
If Condition \emph{$\bold{C4)}$} holds then using Lemma
$\ref{lemma-borne-de-K}$ there exists $\gamma>0$ such that for any
$x\in[0,1]^{d}$ we have $E(S_{n}(x)^{2})\leq
\gamma\prod_{k=1}^{d}[n(x_{k}+h_{n})]$. Let $x\in [0,1]^{d}$ be fixed, using
Lemma $\ref{lemma-borne-de-K}$, there exists $c>0$ such that
\begin{align*}
\|V_{n}(x)\|_{2}&=\frac{\|S_{n}(x)\|_{2}}{\sum_{i\in\Lambda_{n}}a_{i}(x)}\\
&\leq\frac{\sqrt{\gamma}}{c}\left(\prod_{k=1}^{d}[n(x_{k}+h_{n})]\right)^{-1/2}\\
&\leq\frac{\sqrt{\gamma}}{c\,([nh_{n}])^{d/2}}\\
&\leq\frac{2^{d/2}\sqrt{\gamma}}{c\,(nh_{n})^{d/2}}\qquad\textrm{for $n$ sufficiently large.}
\end{align*}
The proof of Theorem $\ref{rate-variance-proba}$ is
complete.
\subsection{Proof of Proposition $\textbf{\ref{rate-biais}}$}
Since $g\in\,\textrm{Lip}(B)$, it follows that
\begin{align*}
\vert B_{n}(x)\vert
&=\bigg\vert\frac{\sum_{i\in\Lambda_{n}}(g(i/n)-g(x))a_{i}(x)}{\sum_{i\in\Lambda_{n}}a_{i}(x)}\bigg\vert\\
&\leq
Bh_{n}\frac{\sum_{i\in\Lambda_{n}}\|(i/n-x)/h_{n}\|a_{i}(x)}{\sum_{i\in\Lambda_{n}}a_{i}(x)}\\
&\leq Bh_{n}.
\end{align*}
The proof of Proposition $\ref{rate-biais}$ is
complete.
\subsection{Proof of Corollary $\textbf{\ref{cor-optimal-rate-ps}}$}
Let $h_{n}=(n^{-d}\log n)^{1/(2+d)}$ then Proposition
$\ref{rate-biais}$ gives
\begin{equation}\label{optimal-rate-biais}
\sup_{x\in [0,1]^{d}}\sup_{g\in\,\textrm{Lip}(B)}\vert
Eg_{n}(x)-g(x)\vert= O\left[\left(\frac{\log
n}{n^d}\right)^{\frac{1}{2+d}}\right].
\end{equation}
Assume that \emph{$\bold{C1)}$} holds. Noting that
$$
\frac{(\log n)^{1/2}}{(nh_{n})^{d/2}}=\left(\frac{\log
n}{n^d}\right)^{\frac{1}{2+d}}
$$
and using ($\ref{rvps3}$) we obtain
\begin{equation}\label{optimal-rvps3}
\sup_{x\in[0,1]^{d}}\vert
g_{n}(x)-Eg_{n}(x)\vert=O_{a.s.}\left[\left(\frac{\log
n}{n^d}\right)^{\frac{1}{2+d}}\right].
\end{equation}
Combining ($\ref{optimal-rate-biais}$) and ($\ref{optimal-rvps3}$)
we derive ($\ref{optimal-rps3}$).\\
\\
Assume that \emph{$\bold{C2)}$} holds for some $0<q<2$. Noting
that
$$
\frac{(\log n)^{1/q}}{(nh_{n})^{d/2}}=\left(\frac{\log
n}{n^d}\right)^{\frac{1}{2+d}}\times(\log n)^{(2-q)/2q}
$$
and using ($\ref{rvps2}$) we obtain
\begin{equation}\label{optimal-rvps2}
\sup_{x\in[0,1]^{d}}\vert
g_{n}(x)-Eg_{n}(x)\vert=O_{a.s.}\left[\left(\frac{\log
n}{n^d}\right)^{\frac{1}{2+d}}\times(\log n)^{(2-q)/2q}\right].
\end{equation}
Combining ($\ref{optimal-rate-biais}$) and ($\ref{optimal-rvps2}$)
we derive ($\ref{optimal-rps2}$).\\
\\
Let $\varepsilon>0$ be fixed and assume that \emph{$\bold{C3)}$}
holds for some $p>2$ which satisfies condition ($\ref{min-p}$).
Applying the point $\bold{3)}$ of Theorem $\ref{rate-variance-ps}$
with $\theta_{1}=1/(2+d)$ and $\theta_{2}=d/(2+d)$ and noting that
$$
v_{n}=\frac{n^{a}(\log
n)^{b}}{(nh_{n})^{d/2}}=n^{\varepsilon}\left(\frac{\log
n}{n^d}\right)^{\frac{1}{2+d}}\Longleftrightarrow\left\{a=\varepsilon\,\,\textrm{and}\,\,b=\frac{1}{2}\right\}
$$
it follows
\begin{equation}\label{optimal-rvps1}
\sup_{x\in[0,1]^{d}}\vert
g_{n}(x)-Eg_{n}(x)\vert=O_{a.s.}\left[n^\varepsilon\left(\frac{\log
n}{n^d}\right)^{\frac{1}{2+d}}\right].
\end{equation}
Combining ($\ref{optimal-rate-biais}$) and ($\ref{optimal-rvps1}$)
we derive ($\ref{optimal-rps1}$). The proof of Corollary
$\ref{cor-optimal-rate-ps}$ is complete.
\subsection{Proof of Corollary $\textbf{\ref{cor-optimal-rate-proba}}$}
Let $h_{n}=n^{-d/(2+d)}$ then Proposition $\ref{rate-biais}$ gives
\begin{equation}\label{optimal-rate-biais-bis}
\sup_{x\in [0,1]^{d}}\sup_{g\in\,\textrm{Lip}(B)}\vert
Eg_{n}(x)-g(x)\vert= O\left[n^{-\frac{d}{2+d}}\right].
\end{equation}
Let $\varepsilon>0$ be fixed and assume that \emph{$\bold{C3)}$}
holds for some $p>2$ which satisfies condition
($\ref{min-p-bis}$). Applying the point $\bold{1)}$ of Theorem
$\ref{rate-variance-proba}$ and noting that
$$
v_{n}=\frac{n^{a}}{(nh_{n})^{d/2}}=n^{-\frac{d}{2+d}+\varepsilon}
\Longleftrightarrow a=\varepsilon
$$
it follows that
\begin{equation}\label{optimal-rvproba2}
\left\|\sup_{x\in[0,1]^{d}}\vert
g_{n}(x)-Eg_{n}(x)\vert\right\|_{p}=O\left[n^{-\frac{d}{2+d}+\varepsilon}\right].
\end{equation}
Combining ($\ref{optimal-rate-biais-bis}$) and
($\ref{optimal-rvproba2}$) we derive ($\ref{optimal-rproba2}$).\\
\\
Since $h_{n}=n^{-d/(2+d)}$ then $(nh_{n})^{-d/2}=h_{n}$. So, if
\emph{$\bold{C4)}$} holds then combining
($\ref{optimal-rate-biais-bis}$) and ($\ref{rvL2}$) we derive
($\ref{optimal-rproba1}$). The proof of Corollary
$\ref{cor-optimal-rate-proba}$ is complete.

\subsection{Proof of Corollary $\textbf{\ref{cor-optimal-rate-ps-mixing}}$}
Let $p>2$ be fixed. Using Rio's inequality \cite{Rio} (see also
Dedecker \cite{JD-tflc}) we obtain the bound
\begin{equation}\label{dedecker-rio-inequality}
\|\varepsilon_{k}E_{\vert
k\vert}(\varepsilon_{0})\|_{\frac{p}{2}}\leq
4\left(\int_{0}^{\alpha_{1,\infty}(\vert
k\vert)}Q_{\varepsilon_{0}}^{p}(u)\,du\right)^{2/p}
\end{equation}
hence condition \emph{$\bold{C^{'}3)}$} is more restrictive
than condition \emph{$\bold{C3)}$}.\\
By Serfling's inequality (see McLeish \cite{McLeish} or Serfling
\cite{Serf}) we know that
$$
\|\varepsilon_{k}E_{\vert k\vert}(\varepsilon_{0})\|_{\infty}\leq
2\|\varepsilon_{0}\|_{\infty}^{2}\phi_{\infty,1}(\vert k\vert)
$$
so condition \emph{$\bold{C^{'}1)}$} is more restrictive than
condition \emph{$\bold{C1)}$}.\\
Now for $0<q<2$ there exists $C(q)>0$ (cf. Inequality (17) in
\cite{MEM-KK}) such that
\begin{equation}\label{majoration-mixing1}
\left\|\sqrt{\vert \varepsilon_{k}E_{\vert
k\vert}(\varepsilon_{0})\vert}\right\|^{2}_{\psi_{\beta(q)}}\leq
C(q)\sqrt{\phi_{\infty,1}(\vert k\vert)}.
\end{equation}
In \cite{MEM-KK} we used the following lemma which can be obtain
by the expansion of the exponential function.
\begin{lemma}\label{lemma-va}
Let $\beta$ be a positive real number and $Z$ be a real random
variable. There exist positive universal constants $A_{\beta}$ and
$B_{\beta}$ depending only on $\beta$ such that
$$
A_{\beta}\,\sup_{p>2}\frac{\| Z\|_{p}}{p^{1/\beta}}\leq\|
Z\|_{\psi_{\beta}}\leq B_{\beta}\,\sup_{p>2}\frac{\|
Z\|_{p}}{p^{1/\beta}}.
$$
\end{lemma}
Consider the coefficient $c_{k}(\beta)$ given by
($\ref{definition-ck}$) and denote
$$
d_{k}(p)=\left(\int_{0}^{\alpha_{1,\infty}(\vert
k\vert)}Q_{\varepsilon_{0}}^{p}(u)\,du\right)^{1/p}
$$
then the following version of lemma $\ref{lemma-va}$ holds.
\begin{lemma}\label{lemma-nonva}
Let $\beta$ be a positive real number. There exist positive
universal constants $A_{\beta}$ and $B_{\beta}$ depending only on
$\beta$ such that for any $k\in\Z^{d}$
$$
A_{\beta}\sup_{p>2}\frac{d_{k}(p)}{p^{1/\beta}}\leq
c_{k}(\beta)\leq B_{\beta}\sup_{p>2}\frac{d_{k}(p)}{p^{1/\beta}}.
$$
\end{lemma}
Now combining lemmas $\ref{lemma-va}$ and $\ref{lemma-nonva}$ and
inequality ($\ref{dedecker-rio-inequality}$) there exists
$C^{'}(q)>0$ such that
\begin{equation}\label{majoration-mixing2}
\left\|\sqrt{\vert \varepsilon_{k}E_{\vert
k\vert}(\varepsilon_{0})\vert}\right\|^{2}_{\psi_{\beta(q)}}\leq
C^{'}(q)\,c_{k}^{2}(\beta(q)).
\end{equation}
Finally condition \emph{$\bold{C^{'}2)}$} is more restrictive
than condition \emph{$\bold{C2)}$} and the proof of Corollary $\ref{cor-optimal-rate-ps-mixing}$ is complete.\\
\\
\textbf{Aknowledgements}. I would like to express my thanks to the anonymous referee for his/her careful
reading of the manuscript and valuable suggestions. I am indebted for \'E. Youndje for many stimulating
conversations on nonparametric estimation.
\bibliographystyle{plain}
\bibliography{xbib}
Mohamed EL MACHKOURI\\
Laboratoire de Math\'ematiques Rapha\"el Salem\\
UMR 6085, Universit\'e de Rouen\\
Site Colbert,\\
F76821 Mont-Saint-Aignan Cedex\\
email : mohamed.elmachkouri@univ-rouen.fr
\end{document}